\documentclass{ifacconf}

\usepackage{graphicx}      
\usepackage{natbib}        
\usepackage{subfigure}
\usepackage{amsmath}
\usepackage{amsfonts}
\usepackage{xcolor}
\usepackage[normalem]{ulem}

\newcommand{\norm}[1]{\left\lVert#1\right\rVert}

\newcommand{\Ltv}{L^{2}(0,T;L^{\infty}(\Omega)^2)}
\newcommand{\Ltk}{L^{2}(0,T;L^{\infty}(\Omega))}

\begin{document}
\begin{frontmatter}

\title{Indirect Optimal Control of Advection-Diffusion Fields through \\ 
Robotic Swarms\thanksref{footnoteinfo}} 

\thanks[footnoteinfo]{This work was supported by the Italian Ministry of Education, University and Research (MIUR).}

\author[First]{Carlo Sinigaglia} 
\author[Second]{Andrea Manzoni} 
\author[Third]{Francesco Braghin} 
\author[Fourth]{Spring Berman}

\address[First]{Department of Mechanical Engineering, Politecnico di Milano, 
   Milano, Italy (e-mail: carlo.sinigaglia@polimi.it).}
\address[Second]{MOX Department of Mathematics, Politecnico di Milano, 
Milano, Italy (e-mail: andrea1.manzoni@polimi.it).}
\address[Third]{Department of Mechanical Engineering, Politecnico di Milano, 
   Milano, Italy (e-mail: francesco.braghin@polimi.it).}
\address[Fourth]{School for Engineering of Matter, Transport and Energy, Arizona State University, Tempe, Arizona 
(e-mail: spring.berman@asu.edu).}

\begin{abstract}                
In this paper, we consider the problem of optimally guiding a large-scale swarm of underwater vehicles that is tasked with the indirect control of an advection-diffusion environmental 
field. The microscopic vehicle dynamics are 
governed by a stochastic differential equation with drift. The drift terms model the self-propelled velocity
of the vehicle 
and 
the velocity field of the currents. In the mean-field setting, the macroscopic vehicle dynamics are 
governed by a Kolmogorov forward equation in the form of a linear parabolic advection-diffusion equation.
The environmental field is governed by an advection-diffusion equation in which 
the advection term is defined by the 
fluid velocity field. The vehicles are equipped with on-board actuators that enable the swarm 
to act as a distributed source in the environmental field, modulated by a scalar control parameter that determines the local source intensity. In this setting, we formulate an optimal control problem to compute the vehicle velocity and actuator intensity fields that drive the environmental field to a desired distribution within a specified amount of time.  
In other words, we design optimal vector and scalar actuation fields to indirectly control the environmental field through a distributed source, produced by the 
swarm. After proving an existence result for the solution of the optimal control problem, we discretize and solve the problem using the Finite Element Method (FEM). The FEM discretization naturally provides an operator that represents the bilinear way in which the controls enter into the dynamics of the vehicle swarm and the environmental field. 
 Finally, we show through numerical simulations the effectiveness of our control strategy in regulating the environmental field to zero or to a desired distribution in the presence of a double-gyre flow field.
\end{abstract}


\end{frontmatter}

\section{Introduction}
Given 
the exponential decrease in cost 
of electronic components over the last few decades, large collectives of robots, or {\it robotic swarms},   
are becoming a viable option for a variety of missions of ever-increasing complexity, such as coverage, mapping, search-and-rescue, and surveillance 
\citep{Dorigo2021}. Controllers for 
robotic swarms should 
satisfy mission requirements while scaling gracefully as the number of robots $N$ 
increases. Mean-field models of robotic swarms (see, e.g., \cite{Elamvazhuthi2020}) describe a swarm as a set of probability densities over space and time; since these models are independent of $N$, they can be used to design controllers for arbitrarily large swarms, with the caveat that the  distribution of the swarm is controlled rather than individual robots.  

In this paper, we take advantage of the mean-field model's invariance to swarm size by using such a model to design scalable robotic swarm controllers that achieve indirect control of a distributed process that evolves according to a 
Partial Differential Equation (PDE), such as the concentration field of a contaminant in a fluid flow. The robots are underwater vehicles that are each equipped with an actuator that acts as a source for the process, and we aim to indirectly control the process through the coordinated motion of the swarm and the source actuation. A similar problem has previously been considered for finite teams of mobile robots, in which individual robot trajectories are controlled.   
For example, \cite{ChengPaley2021} 
present an optimal control approach that uses an operator-valued Riccati equation to formulate the optimal actuation as a function of the optimal guidance, and then recast the problem in terms of the latter  
alone to jointly optimize the guidance of the 
robots and their associated actuation.
\cite{demetriou2021controlling} describes a path-dependent reachability approach that takes into account constraints on the robots' motion and real-time implementation while regulating a spatially distributed process using local decentralized measurements only. 
 
We define the indirect control problem for a distributed robotic swarm whose mean-field dynamics
are 
modeled by a Kolmogorov forward equation in the form of a linear parabolic advection-diffusion PDE.
We formulate an optimal control problem (OCP) for this mean-field model, which is coupled with the advection-diffusion dynamics of the environmental field and prove an existence theorem for the OCP  using techniques for optimal control of PDEs.
Then, we derive a set of first-order necessary optimality conditions and solve them 
numerically using a Finite Element Method (FEM) discretization. Finally, we evaluate the effectiveness of our control strategy through numerical simulations of regulation and target tracking problems 
in the presence of a double-gyre flow field.

\section{Problem formulation}
\begin{figure}
    \centering
    \includegraphics[width=0.65\linewidth]{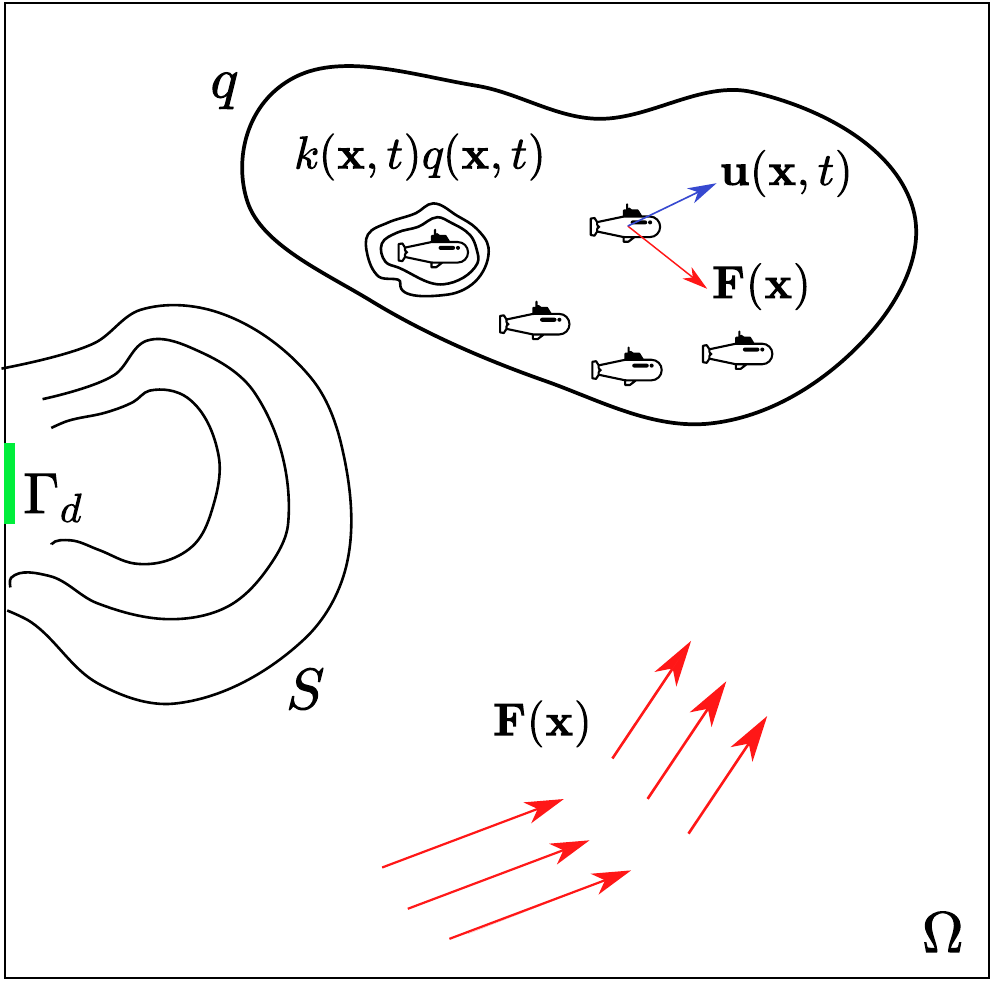}
    \caption{Illustration of the indirect control problem. The environmental field $S$ diffuses from the segment of the boundary  $\Gamma_d$ under the advection of a fluid 
    velocity field $\mathbf{F}$. The vehicle velocity field 
    $\mathbf{u}$ and 
    actuator intensity field 
    $k$ are designed to drive $S$ to a target distribution.
    }
    \label{schem_diag} 
\end{figure}

We consider a swarm of 
robots, labeled $i = 1,...,N$, that move in a bounded fluid domain $\Omega \in \mathbb{R}^2$. Robot $i$ occupies position $\mathbf{X}_i(t) \in \Omega$ at time $t$ and moves with velocity $\mathbf{v}(\mathbf{x},t) \in \mathbb{R}^2$, 
which is the sum of its self-propelled velocity $\mathbf{u}(\mathbf{x},t)$ and the fluid velocity at its position,  $\mathbf{F}(\mathbf{x})$. This motion is perturbed by a 2-dimensional Wiener process $\mathbf{W}(t)$, which models stochasticity arising from inherent sensor and actuator noise or intentionally programmed ``diffusive'' exploratory behaviors. The robot's position evolves according to
the following Stochastic Differential Equation:
\begin{equation*}
\begin{cases}
d \mathbf{X}_i(t) & = ~\mathbf{v}(\mathbf{X}_i,t)dt +\sqrt{2 D_q} d \mathbf{W}(t)+\mathbf{n}(\mathbf{X}_i(t)) d \psi(t) \\ \mathbf{X}_i(0) & = ~\mathbf{X}_{i,0},
\end{cases}
\end{equation*}
where 
$D_q>0$ is a diffusion coefficient, 
$\mathbf{n}(\mathbf{x})$ is the unit normal to the domain boundary $\partial \Omega$ at $\mathbf{x}$, 
and $\psi(t) \in \mathbb{R}$ is a reflecting function, 
which ensures that the 
swarm does not exit the domain.
We assume that each robot $i$ carries an on-board actuator that acts as a source of intensity $k_i(t) \in \mathbb{R}$ in a scalar environmental field $S(\mathbf{x},t)$, $\mathbf{x} \in \Omega$. We also assume that part of the boundary, $\Gamma_d \in \partial\Omega$, acts as a source 
with constant intensity $S_d$ for simplicity. The field $S$ is an 
advection-diffusion process with diffusion coefficient $D_S >0$, modeled by the following PDE problem: 
\begin{equation*} 
\begin{array}{ll}
\displaystyle \frac{\partial S}{\partial t} - D_S \Delta S + \mathbf{F}\cdot \nabla S= \sum_{i=1}^{N} k_i(t) \delta(\mathbf{x}-\mathbf{X}_i(t))  &  \textrm{in} \quad \Omega \\
\vspace{-3mm} \\
S = S_d \,  1_{\Gamma_d}  \quad\quad\quad \textrm{on} \quad \partial \Omega & \\
S(\mathbf{x},0) = S_0(\mathbf{x})  \quad \textrm{on} \quad \Omega \times \{ 0 \}, & \\
\end{array}
\end{equation*}
where $1_{\Gamma_d}$ is an indicator function. The right-hand side of the PDE consists of 
the cumulative effect of the point 
sources, which have the same dynamics as the robots since they are installed on-board.
Note 
that the explicit dependence of $S$ on space and time is omitted when clear from the context. 


In the limit as $N \to \infty$, we obtain a mean-field model that describes the evolution of the probability density $q(\mathbf{x},t)$ of a single robot occupying position $\mathbf{x}$ at time $t$, or alternatively, the swarm density at this position and time. For the robot dynamics we consider here, this model takes the form of 
a linear parabolic advection-diffusion problem with no-flux boundary conditions; see, e.g., \cite{sinigaglia2021density}. The coupled system dynamics are therefore:  
\begin{equation} 
\label{dyn}
\begin{array}{ll}
\displaystyle \frac{\partial q}{\partial t} + \nabla \cdot (-D_q \nabla q + \mathbf{u}q + \mathbf{F}q) = 0 &  \textrm{in} \quad \Omega \\
\vspace{-3mm} \\
(-D_q \nabla q + \mathbf{u}q + \mathbf{F}q )\cdot \mathbf{n} = 0 & \textrm{on} \quad \partial \Omega \\
q(\mathbf{x},0) = q_0(\mathbf{x}) & \textrm{in} \quad \Omega \times \{ 0 \}, \\
\vspace{-1mm} \\
\displaystyle \frac{\partial S}{\partial t} - D_S \Delta S + \mathbf{F}\cdot \nabla S= k q  &  \textrm{in} \quad \Omega \\
\vspace{-3mm} \\
S = S_d \,  1_{\Gamma_d} &  \textrm{on} \quad \partial \Omega \\
S(\mathbf{x},0) = S_0(\mathbf{x}) & \textrm{on} \quad \Omega \times \{ 0 \}. \\
\end{array}
\end{equation}
Note that the right-hand side of the $S$ dynamics is now the product of the intensity of the distributed actuation $k$ and the local swarm density $q$. As a consequence of this Eulerian perspective, $k = k(\mathbf{x},t)$ varies in both space and time. The robotic 
swarm thus aims at controlling the concentration $S$ of an advection-diffusion field through a localized actuation field $kq$ that depends on both the actuator intensity $k$ and the swarm density $q$ at each spatial location. See Figure \ref{schem_diag} for an illustration of the 
problem.

Defining the state variables as $q$ and $S$,
the state dynamics consists of the one-way coupled system of PDEs \eqref{dyn}. The objective of the control problem is then to find the optimal actuation for $\mathbf{u}(\mathbf{x},t)$ and $k(\mathbf{x},t)$ to guide the space-time evolution of the field $S$ to a target distribution $z$, which may be static or dynamic, at a given time $T$. This control objective
can be easily 
encoded in a standard quadratic cost functional of the form:
\begin{equation}
\label{cost}
\begin{aligned}
    &J = \frac{\alpha_T}{2}\int_{\Omega} (S(\mathbf{x},T)-z)^2 d\Omega  +\frac{\alpha}{2} \int_{0}^{T}\int_{\Omega} (S-z)^2 d\Omega d t \\
    & ~~ + \frac{\beta}{2} \int_{0}^{T}\int_{\Omega} \norm{\mathbf{u}}^2 d\Omega d t + \frac{\gamma}{2} \int_{0}^{T}\int_{\Omega} k^2 d\Omega d t,
\end{aligned}
\end{equation}
where 
$\alpha_T,\alpha,\beta,\gamma \geq 0$ are weighting constants. In the regulation problem, for example, we set $z=0$ and 
seek 
the optimal trade-off between effectively regulating the field $S$ and 
minimizing the overall energy expended by the swarm for propulsion motion and source actuation.  

\section{The Optimal Control Problem}
In this section, we prove the existence of optimal controls, 
derive a system of first-order necessary optimality conditions using the Lagrangian method, and 
provide a consistent discretization of the OCP using the FEM.

\subsection{Analysis}
Both the $q$ and $S$ dynamics satisfy rather standard advection-diffusion equations of linear parabolic type; see, e.g., \cite[Chapter 7]{MQS}. 
The natural functional space for the swarm density function $q$ which is subjected to zero-flux Neumann boundary conditions is $q \in L^2(0,T,H^1(\Omega))$, while the Dirichlet boundary suggests the choice of  $\mathring{S} \in L^2(0,T,H^1_0(\Omega))$ for the ``lifted'' field variable $\mathring{S} = S - \tilde{S}_d$, where $\tilde{S}_d$ is a suitable extension of the boundary datum to the domain $\Omega$ -- see, e.g., \cite[Chapter 8]{salsa2015partial}. 
It is also standard to select $\frac{\partial q}{\partial t} \in L^2(0,T,H^1(\Omega)^{*})$ and $\frac{\partial \mathring{S}}{\partial t} \in L^2(0,T,H_0^1(\Omega)^{*})$ so that the functional space for $q$ is actually $ H^1(0,T;H^1(\Omega),H^1(\Omega)^{*}) = \{ y \in L^2(0,T;H^1(\Omega)) : \dot{y} \in L^2(0,T;H^1(\Omega)^{*}\}$, and the same holds true for $\mathring{S}$, substituting $H^1$ with $H_0^1$. Therefore, we set 
$\mathcal{Y} = H^1(0,T;H^1(\Omega),H^1(\Omega)^{*}) \times H^1(0,T;H_0^1(\Omega),H_0^1(\Omega)^{*})$ as the state space, i.e., $y=(q,S) \in \mathcal{Y}$. As done in, e.g., \cite{roy2018fokker} and \cite{sinigaglia2021density} for similar problems involving the Kolmogorov forward equation alone, we consider $L^{\infty}$ spaces for the control fields for which energy-like inequalities are readily available; that is, we select $\mathcal{U} = \Ltv \times \Ltk$ as the control space, so that $v = (\mathbf{u},k) \in \mathcal{U}$.
Besides the choice of the functional spaces for states and controls, we make the following standard assumptions: 
\begin{subequations}
\begin{align}
    \mathbf{F} \in L^{\infty}(\Omega)^2 \tag{A1} \label{a1} \\
    D_S,D_q>0 \tag{A2} \label{a2} \\
    S_d\, \, \text{is bounded} \tag{A3}\label{a3} \\
    q_0,S_0 \in L^2(\Omega) \tag{A4} \label{a4}
\end{align}
\end{subequations}

The weak formulation of the  
PDE problem governing the swarm dynamics is: find $q \in L^2(0,T;H^1(\Omega))$ such that for a.e. $t \in (0,T)$,
\begin{equation}
\begin{aligned}
&\int_{\Omega} \frac{\partial q}{\partial t} \phi d \,\Omega + a_q(q,v;\mathbf{u}) = 0
\vspace{1mm} \\
&q(0)=q_{0} &
\end{aligned}
\end{equation}
for every $\phi \in H^1(\Omega)$, where
\begin{equation*}
    a_q(q,\phi;\mathbf{u}) = \int_{\Omega}(D_q \nabla q \cdot \nabla \phi-(\mathbf{u}+\mathbf{F}) \cdot \nabla \phi q) d \Omega.
\end{equation*}
 The weak formulation of the 
 PDE problem for the ``lifted'' variable $\mathring{S}$ is: 
 find $\mathring{S} \,\in L^2(0,T;H_0^1(\Omega))$ such that for a.e. $t \in (0,T)$, 
\begin{equation}
\begin{aligned}
&\int_{\Omega} \frac{\partial \mathring{S}}{\partial t} \phi \,d \Omega +a_S(\mathring{S},\phi)=\int_{\Omega} \, k q\,\phi \,d\Omega -a_S(\tilde{S}_d,\phi)
\vspace{1mm} \\
&\mathring{S}(0)=S_{0} & 
\end{aligned}
\end{equation}
for every $\phi \in H_0^1(\Omega)$, where 
\begin{equation*}
    a_S(S,\phi) = \int_{\Omega} D_S \nabla S \cdot \nabla \phi + \mathbf{F} \cdot \nabla S \,\phi d\Omega.
\end{equation*}
We also define the linear functional $F \in H^{-1}=H_0^1(\Omega)^{*}$ as $F\phi = \int_{\Omega} k\,q\,\phi\,d\Omega - a_S(\tilde{S}_d,\phi)$. In the following, we will need a bound on the operator norm of $F$, which we prove in the lemma below.
\begin{lem}[Bound on $F$]
\label{F_bound} Let assumptions \eqref{a1}, \eqref{a2}, \eqref{a3}, and \eqref{a4} hold. Then the following bound on the norm of $F$ holds:
\begin{equation*}
\hspace{-4mm} \norm{F}_{H^{-1}(\Omega)} \leq C_p \norm{kq}_{L^2(\Omega)} + (C_p\norm{\mathbf{F}}_{L^{\infty}(\Omega)^2}+D_S)\norm{\nabla \tilde{S}_d}_{L^2(\Omega)},
\end{equation*}
where $C_p>0$ is the Poincarè inequality constant.
\begin{pf}
From the definition of $F$ and the Cauchy-Schwarz and Poincar\'e inequalities, we have: 
\begin{equation*}
\begin{aligned}
    &|F\phi| \leq \Big( \norm{kq}_{L^2(\Omega)}+\norm{\mathbf{F}}_{L^{\infty}(\Omega)^2}\norm{\nabla \tilde{S}_d}_{L^2(\Omega)}\Big) \norm{\phi}_{L^2(\Omega)} \\
    & + D_S \norm{\nabla \tilde{S}_d}_{L^2(\Omega)} \norm{\nabla \phi}_{L^2(\Omega)} \\ 
    & \leq C_p \Big( \norm{kq}_{L^2(\Omega)}+\norm{\mathbf{F}}_{L^{\infty}(\Omega)^2}\norm{\nabla \tilde{S}_d}_{L^2(\Omega)}\Big)\norm{\nabla \phi}_{L^2(\Omega)}\\
    &+ D_S \norm{\nabla \tilde{S}_d}_{L^2(\Omega)}\norm{\nabla \phi}_{L^2(\Omega)}.
\end{aligned}
\end{equation*}
Regrouping the terms and using the definition of operator norm in $H^{-1}(\Omega)$, the result follows. \qed
\end{pf}
\end{lem}

Existence and well-posedness of the state dynamics follow from the well-posedness of the $q$ dynamics and basic energy estimates on the $S$ dynamics. This is a consequence of the one-way coupling from $q$ to $S$. Following the same arguments as in \cite{sinigaglia2021density}, we have that 
\begin{equation*}
    \norm{q}^2_{H^1(0,T;H^1(\Omega),H^1(\Omega)^{*})} \leq C_0(\norm{\mathbf{u}}^2)\norm{q_0}^2_{L^2(\Omega)}.
\end{equation*}
To prove the well-posedness of the $S$ dynamics, we note that $kq \in L^2(0,T;L^2(\Omega))$ since
\begin{equation*}
\begin{aligned}
    &\norm{kq}_{L^2(0,T;L^2(\Omega))}^2 \leq \norm{k}_{L^2(0,T;L^{\infty}(\Omega))}^2\norm{q}_{L^2(0,T;L^{2}(\Omega))}^2 \\
    & \hspace{5mm} \leq C \norm{k}_{L^2(0,T;L^{\infty}(\Omega))}^2\norm{\mathbf{u}}_{L^2(0,T;L^{\infty}(\Omega)^2)}^2,
\end{aligned}
\end{equation*}
and the latter quantity in the inequality is bounded by the definition of the control space $\mathcal{U}$. Therefore, $S$ satisfies an advection-diffusion equation with $L^2$ right-hand side for which existence and uniqueness results are available -- see, e.g., \cite[Theorem 9.9]{salsa2015partial}. 

A number of standard a priori estimates can be derived for the $S$ dynamics as well; see,  e.g., \cite[Theorem 7.1]{MQS}. In particular, it can be shown that
\begin{equation*}
    \hspace{-4.5mm} \norm{S}_{L^2(0,T;H_0^1(\Omega))}^2 \leq \frac{e^{2\lambda T}}{\alpha_0}\Big(\norm{S_0}_{L^2(\Omega)}^2 + \frac{1}{\alpha_0}\norm{F}_{L^2(0,T;H^{-1}(\Omega))}^2 \Big),
\end{equation*}
where $\lambda = \frac{\norm{\mathbf{F}}^2_{L^{\infty}(\Omega)^2}}{D_S}$ and $\alpha_0 = \min\{ \frac{D_S}{2},\frac{\norm{\mathbf{F}}^2_{L^{\infty}(\Omega)^2}}{2\,D_S} \}$. From Lemma \ref{F_bound} and the bounds on $\norm{kq}_{L^2(0,T;L^2(\Omega)}$, it is 
clear that $\norm{S}_{L^2(0,T;H_0^1(\Omega))}$ is bounded by the control norms on $\mathbf{u}$ and $k$. Regarding $\dot{S}=\frac{\partial S}{\partial t}$, we have that
\begin{eqnarray*}
   \norm{\dot{S}}_{L^2(0,T;H^{-1}(\Omega))}^2   ~\leq & \hspace{-12mm} C_0 \norm{S_0}_{L^2(\Omega)}^2 + \\ 
   & ~~~~ (\frac{C_0}{\alpha_0}+2)\norm{F}_{L^2(0,T;H^{-1}(\Omega))}^2, 
\end{eqnarray*}
where $C_0 = \frac{2}{\alpha_0}e^{2\lambda T}(D_S+
C_p\norm{\mathbf{F}}_{L^{\infty}(\Omega)^2})^2$,
since both $\norm{\dot{S}}_{L^2(0,T;H^{-1}(\Omega))}$ and $\norm{S}_{L^2(0,T;H_0^1(\Omega))}$ are bounded by the control norms. Therefore, we can conclude that 
\begin{equation*}
\begin{aligned}
    &\norm{S}^2_{H^1(0,T;H_0^1(\Omega),H^{-1}(\Omega))} \\
    &=   \norm{\dot{S}}_{L^2(0,T;H^{-1}(\Omega))}^2 + \norm{S}^{2}_{L^2(0,T;H_0^1(\Omega))} \\
    &\leq f(\norm{k}_{L^2(0,T;L^{\infty}(\Omega))},\norm{\mathbf{u}}_{L^2(0,T;L^{\infty}(\Omega)^2)}),
    \end{aligned}
\end{equation*} 
which will turn out to be useful in the proof of existence of optimal controls.

We define the control-to-state operator as the map $(S,q) = \Xi[\mathbf{u},k]$ which associates to each control function $v \in \mathcal{U}$ with a  corresponding state $y \in \mathcal{Y}$. The following result regarding the Fr\'echet differentiability of the control-to-state operator is also needed 
to prove existence of optimal controls.

\begin{lem}[Differentiability of the control-to-state map]
Let assumptions \eqref{a1}, \eqref{a2}, \eqref{a3}, and \eqref{a4} hold. Then the control-to-state map $(S,q) = \Xi[\mathbf{u},k]$ is Fr\'echet differentiable and the directional derivative $(z_S,z_q) =\Xi'[\mathbf{u},k](\mathbf{h},l)$ at $(\mathbf{u},k) \in \mathcal{U}$ in the direction $(\mathbf{h},l) \in \mathcal{U}$ is the solution of the coupled PDE system:
\begin{equation*}
\begin{array}{ll}
\displaystyle \frac{\partial z_q}{\partial t} + \nabla \cdot (-D_q \nabla z_q + \mathbf{u}z_q + \mathbf{F}z_q) = -\nabla \cdot (\mathbf{h}\,q) &  \textrm{in} \quad \Omega \\
\vspace{-3mm} \\
(-D_q \nabla z_q + \mathbf{u}z_q + \mathbf{F}z_q )\cdot \mathbf{n} = -\mathbf{h}\cdot \mathbf{n}\,q & \textrm{on} \quad \partial \Omega \\
\vspace{-2mm} \\
z_q(\mathbf{x},0) = 0 & \hspace{-6mm} \textrm{in} \quad \Omega \times \{ 0 \}, \\
\vspace{1mm} \\
\displaystyle \frac{\partial z_S}{\partial t} - D_S \Delta z_S + \mathbf{F}\cdot \nabla z_S= k z_q + l q  &  \textrm{in} \quad \Omega \\
\vspace{-2mm} \\
z_S = 0 \, &  \textrm{on} \quad \partial \Omega \\
z_S(\mathbf{x},0) = 0 & \hspace{-6mm} \textrm{on} \quad \Omega \times \{ 0 \}. \\
\end{array}
\end{equation*}

\begin{pf}(Sketch) The derivation of the equations governing the sensitivity $z_q$ of the swarm dynamics, and bounds on the norm of $z_q$,
can be found e.g. in \cite{roy2018fokker} and \cite{sinigaglia2021density}. On the other hand, the expression for the dynamics of $z_S$ can be obtained by formally computing
the directional derivative, that is, the limit $\lim_{s \to 0}\frac{S(\mathbf{u}+s\mathbf{h},k+sl)-S(\mathbf{u},k)}{s}$. Bounds on $z_S$ are obtained by adapting the above results on the norm of $S$, and noting that by the triangular inequality,
   $ \norm{kz_q+lq}_{L^2(0,T;L^2(\Omega))} \leq \norm{kz_q}_{L^2(0,T;L^2(\Omega))} + \norm{lq}_{L^2(0,T;L^2(\Omega))},$
so that continuity of the sensitivity  $z_q$ and $z_S$ with respect to variations of the control functions $\mathbf{h}$, $l$ can be easily obtained, thus proving the differentiability of the control-to-state map. \qed
\end{pf}
\end{lem}

We are now ready to prove a result concerning the existence of optimal controls.  

\begin{thm}[Existence of optimal controls] Let assumptions \eqref{a1}, \eqref{a2}, \eqref{a3}, and \eqref{a4} hold. Then
there exists an optimal control  $\bar{v} = (\bar{\mathbf{u}},\bar{k})$ that minimizes the cost functional \eqref{cost} subject to the dynamics \eqref{dyn}, that is, $\bar{y}=(\bar{S},\bar{q})=\Xi[\bar{\mathbf{u}},\bar{k}]$.
\end{thm}

\begin{pf}    
(Sketch) Existence results for bilinear optimal control problems involving the Kolmogorov forward equation with space-time dependent controls have been proved in \cite{sinigaglia2021density}, \cite{roy2018fokker}, and references therein. Choosing a minimizing sequence $(\mathbf{u}_n,k_n)$, due to the weak*  sequential compactness of the control space, we have that  
\begin{equation*}
\begin{aligned}
    &\mathbf{u}_{n} \stackrel{*}{\rightharpoonup} \bar{\mathbf{u}} \quad \text { (weakly star) in } L^2(0,T;L^{\infty}(\Omega)^2) \\
    &k_{n} \stackrel{*}{\rightharpoonup} \bar{k} \quad \text { (weakly star) in } L^2(0,T;L^{\infty}(\Omega))
\end{aligned}
\end{equation*}
Due to the bounds on $S$ and $q$, we also have that the resulting sequence $(S_n,q_n)=\Xi[\mathbf{u}_n,q_n]$ is bounded and thus weakly convergent in $\mathcal{Y}$ to $(\bar{S},\bar{q}) \in \mathcal{Y}$, see e.g. \cite[Appendix S, Theorem 3]{evans}. It remains to prove that: 
\begin{equation*}
    \int_{0}^{T}\int_{\Omega} \,k_n q_n\,\phi\,d\Omega\,dt \to  \int_{0}^{T}\int_{\Omega} \,\bar{k} \bar{q}\,\phi\,d\Omega\,dt
\end{equation*}
for each $\phi \in L^2(0,T;H_0^1(\Omega))$. To this end, we write 
\begin{equation*}
\begin{aligned}
    &\int_{0}^{T} \int_{\Omega}\left(k_{n} q_{n}-\bar{q} \bar{k}\right) \phi d \Omega d t\\
    &=\int_{0}^{T}\int_{\Omega}  \bar{q} \phi\left(k_{n}-\bar{k}\right) d \Omega d t+\int_{0}^{T} \int_{\Omega}\left(q_{n}-\bar{q}\right) k_{n} \phi \,d \Omega d t.
\end{aligned}
\end{equation*}
Since $\bar{q}\phi \in L^2(0,T;L^1(\Omega))$, 
the dual of $L^2(0,T;L^{\infty}(\Omega))$, and $k_{n} \stackrel{*}{\rightharpoonup} \bar{k}$, 
the first integral tends to zero by Lebesgue's dominated convergence theorem.
To analyze the second integral, 
we use the Aubin-Lions Lemma 
(see, e.g.,  \cite[Appendix A, Theorem A.19]{MQS}) to ensure that  $q_n \to \bar{q}$ strongly in $L^2(0,T;L^2(\Omega))$. Then we obtain:
\begin{equation*}
\begin{aligned}
    &|\int_{0}^{T} \int_{\Omega}\left(q_{n}-\bar{q}\right) k_{n} \phi \,d \Omega d t| \leq \\
    & \hspace{-2mm} \norm{k_n}_{L^2(0,T;L^{\infty}(\Omega))}\norm{\phi}_{L^2(0,T;L^2(\Omega))}\norm{q_n-\bar{q}}_{L^2(0,T;L^2(\Omega))} = 0.
\end{aligned} 
\end{equation*}

Since $\Omega$ is bounded, the weak*  convergence of $\mathbf{u}_n \in L^2(0,T;L^{\infty}(\Omega)^2)$ to some $\bar{\mathbf{u}} \in L^2(0,T;L^{\infty}(\Omega)^2)$ implies weak convergence of $\mathbf{u}_n$ to $\bar{\mathbf{u}}$ in $L^2(0,T;L^{2}(\Omega)^2)$. The same holds for $k$; that is, $k_n$ weakly converges to $\bar{k}$ in $L^2(0,T;L^{2}(\Omega))$. Then, exploiting the fact that $S_n$ weakly converges to $\bar{S}$ in $L^2(0,T;H_0^{1}(\Omega))$ and that $J(S,\mathbf{u},k)$ is convex and continuous in $L^2(0,T,H_0^1(\Omega))\times L^2(0,T;L^{2}(\Omega)^2)\times L^2(0,T;L^{2}(\Omega))$, we conclude that
\begin{equation*}
    J(\bar{S},\bar{\mathbf{u}},\bar{k}) \leq \lim_{n \to \infty} \inf J(S_n,\bar{\mathbf{u}}_n,k_n) = \inf J.
\end{equation*}
Therefore, the pair $(\bar{y},\bar{v})$ is an optimal pair for the considered optimal control problem.   \qed
\end{pf}

We note that uniqueness of an optimal solution is not guaranteed, due to the bilinear way in which both controls $\mathbf{u}$ and $k$ enter into the coupled system dynamics.s

\subsection{Optimality Conditions}

We now derive a system of first-order necessary optimality conditions using the Lagrangian multipliers method. For the problem at hand, the Lagrangian can be defined as
\begin{equation}
\begin{aligned}
    &\mathcal{L} = J - \int_{0}^{T}\int_{\Omega} \lambda_q \left( \frac{\partial q}{\partial t} + \nabla \cdot (-D_q \nabla q + \mathbf{u}q) \right) d\Omega dt \\
    & ~~~~~ + \int_{0}^{T}\int_{\Omega} \lambda_S \left( k q - \frac{\partial S}{\partial t} + D_S \Delta S  \right) d\Omega dt.
\end{aligned}
\end{equation}
Note that we have defined adjoint fields $\lambda_q$ and $\lambda_S$ that are related to the state dynamics of both $q$ and $S$. The adjoint dynamics for $\lambda_q$ and $\lambda_S$ thus satisfy:
\begin{equation*}
\begin{array}{ll}
\displaystyle - \frac{\partial \lambda_q}{\partial t} -D_q \Delta \lambda_q - \mathbf{u} \cdot \nabla \lambda_q = k \, \lambda_S &  \textrm{in} \quad \Omega \\
\vspace{-3mm} \\
\nabla \lambda_q \cdot \mathbf{n} = 0 & \textrm{on} \quad \partial \Omega \\ 
\lambda_q(\mathbf{x},T) = 0 & \textrm{in} \quad \Omega \times \{ T \}, \\
\vspace{-2mm} \\
\displaystyle -\frac{\partial \lambda_S}{\partial t} - D_S \Delta \lambda_S -\mathbf{F} \cdot \nabla \lambda_S = \alpha\,(S-z)  &  \textrm{in} \quad \Omega \\
\vspace{-3mm} \\
\lambda_S = 0 &  \textrm{on} \quad \partial \Omega \\
\lambda_S(\mathbf{x},T) = \alpha_T\,(S(\mathbf{x},T)-z) & \textrm{in} \quad \Omega \times \{ T \}, \\
\end{array}
\end{equation*}
which are obtained by taking the first variation of the Lagrangian along variations in swarm density $q$ and the environmental field $S$, respectively.
Note that the coupling between the adjoint fields is dual with respect to the state dynamics. The coupling is from $q$ to $S$ in the state system, while it is from $\lambda_S$ to $\lambda_q$ in the adjoint system. The dual of the forcing term $k\,q$ in the $S$ dynamics is the forcing term $k \lambda_S$ in the $\lambda_q$ dynamics. 

The reduced gradients of $J$ with respect to $\mathbf{u}$ and $k$ can therefore be expressed as 
\begin{equation}
\begin{array}{l}
    \nabla J_{\mathbf{u}} = \beta   \mathbf{u} + q \, \nabla \lambda_q,  \\
    \nabla J_{k}          = \gamma  k          + q \, \lambda_S, 
\end{array}
\end{equation}
by taking the first variations of the Lagrangian in the directions of $\mathbf{u}$ and $k$, respectively.
Note that, despite $k$ entering linearly into the $S$ dynamics, the reduced gradient $\nabla J_{k}$ depends on the dynamics of the swarm density $q$. This is a consequence of the multiplicative nature of the $S$ forcing term $k\,q$.

\subsection{Numerical Discretization}
The OCP with coupled system dynamics \eqref{dyn} is discretized in the state variables $S$ and $q$ using the Finite Element Method (FEM).
The discretized state dynamics are 
\begin{equation}
    \begin{array}{l}
         M_q \dot{\mathbf{q}} + \Big(A_q - B_{\mathbf{F}}^{\top} - \mathbb{B}^{\top}\mathbf{u}\Big) \mathbf{q} = \mathbf{0}  \\
         M_S \dot{\mathbf{S}} + \Big(A_S-B_{\mathbf{F}}^{\top}\Big) \mathbf{S} = \mathbb{C}\mathbf{k} \mathbf{q} \\
         \mathbf{q}(0)=\mathbf{q}_0 \\
         \mathbf{S}(0)=\mathbf{S}_0,
    \end{array}
\end{equation}
where 
$\mathbb{B}$ is the rank-3 transport coefficient  tensor defined by $\mathbb{B}_{ijk} = \int_{\Omega} \frac{\partial \phi_j}{\partial x}\phi_i\,\phi_k\,d\Omega$; 
$\mathbb{B}\mathbf{u}$ is the tensor vector product defined by $(\mathbb{B}\mathbf{u})_{ij} = \sum_{k=1}^{N_u} \mathbb{B}_{ijk}\,u_k$; 
$\mathbb{C}$ is the reaction tensor defined by $\mathbb{C}_{ijk} = \int_{\Omega}  \phi_j\phi_i\,\phi_k\,d\Omega$; and 
$M$, $B$, and $A$ are the usual FEM mass, transport, and stiffness matrices, respectively. 

The discretization of the adjoint system is: 
\begin{equation}
    \begin{array}{l}
         -M_q \dot{\boldsymbol{\lambda}}_q + \Big(A_q - B_{\mathbf{F}} - \mathbb{B}\mathbf{u} \Big) \boldsymbol{\lambda}_q = \mathbf{k} \mathbb{C}^{\top} \boldsymbol{\lambda}_S   \\
         -M_S \dot{\boldsymbol{\lambda}}_S + A_S \boldsymbol{\lambda}_S= \alpha\,M_S (\mathbf{S}-\mathbf{z}) \\
         \boldsymbol{\lambda}_q(T) = \mathbf{0} \\
         \boldsymbol{\lambda}_S(T) = \alpha_T M_S(\mathbf{S}(T)-\mathbf{z}). \\
    \end{array}
\end{equation}

Finally, the reduced gradient discretization is: 
\begin{equation}
    \begin{array}{l}
    \nabla J_{\mathbf{u}}   = \beta   M_u \mathbf{u} + \mathbf{q}^{\top} \mathbb{B} \boldsymbol{\lambda}_q \\
    \vspace{-4mm} \\
    \nabla J_{\mathbf{k}}   = \gamma  M_k \mathbf{k}   + \mathbf{q}^{\top} \mathbb{C}^{\top} \boldsymbol{\lambda}_S.
\end{array}
\end{equation}

 We can apply the same reasoning as in \cite{sinigaglia2021density} to perform numerical gradient computation. 
 Therefore, we use the Discretize-then-Optimize (DtO) approach (see e.g., \cite[Chapter 8]{MQS}) to numerically solve the problem while avoiding inconsistencies in the gradient computation. In order to do so, the discrete Lagrangian must 
 be computed and differentiated. This computation, which is 
 very similar to the one in our previous work (see \cite{sinigaglia2021density} for a more detailed treatment of 
 the problem for the swarm dynamics alone), is not reported here due to space constraints. 
 
Since the coupling is one-way, at each time step we advance the dynamics of the swarm density  $q$ and then 
solve the problem for $S$. Using similar reasoning, we first solve the discrete adjoint dynamics with respect to $\lambda_S$, and then the adjoint problem for $\lambda_q$. It can be checked that carrying out the optimization at the continuous level and then discretizing the optimality conditions, that is, adopting the Optimize-then-Discretize (OtD) approach, results in the same system of equations at the semi-discrete level; up to choosing a suitable time-discretization, the two approaches fully commute.

\section{Simulation Results}
In this section, we present numerical simulation  results that show the effectiveness of our control strategy. The computational domain $\Omega= [0,1]^2$ is discretized into 
a triangular mesh with $N_s = 2,788$ degrees of freedom, and 
the time interval $[0,T]$ is discretized into  $N_t=61$ time steps, where the final time is $T=1.5(s)$. The resulting fully discrete optimization problem has 510,204 control variables and 170,068 state variables. A steady double-gyre flow field is chosen as the fluid velocity field  $\mathbf{F}$. 

Using the DtO method, the reduced gradient is computed with respect to the control variables only. Computations are 
carried out in MATLAB using  a modified version of the \texttt{redbKit}  (\cite{quarteronireduced}) library to assemble 
the FEM matrices and tensors and the 
\texttt{Tensor Toolbox} (\cite{tensor_toolbox}) 
to perform efficient tensor computations. The nonlinear optimization software \texttt{Ipopt} 
(\cite{wachter2006implementation}) is then used to solve the resulting nonlinear optimization problem.

\begin{figure}[h!]
\centering
\includegraphics[width=0.47\textwidth]{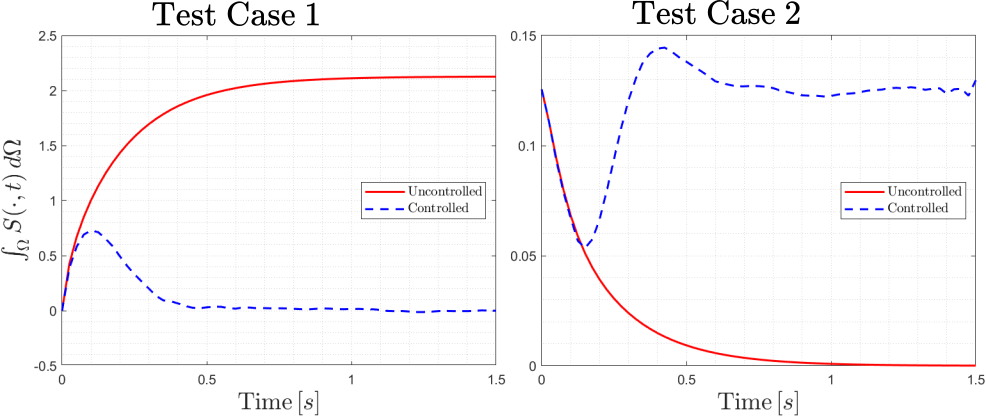}
\caption{ Controlled and uncontrolled dynamics of the total environmental field mass for both test cases.
After a common transient phase during which the robotic swarm moves toward where its actuation is most effective, the controlled dynamics are driven to the total mass of the target distribution $z$. }
\label{mass_conv}
\end{figure}

\begin{figure}[h!]
\centering
\includegraphics[width=0.45\textwidth]{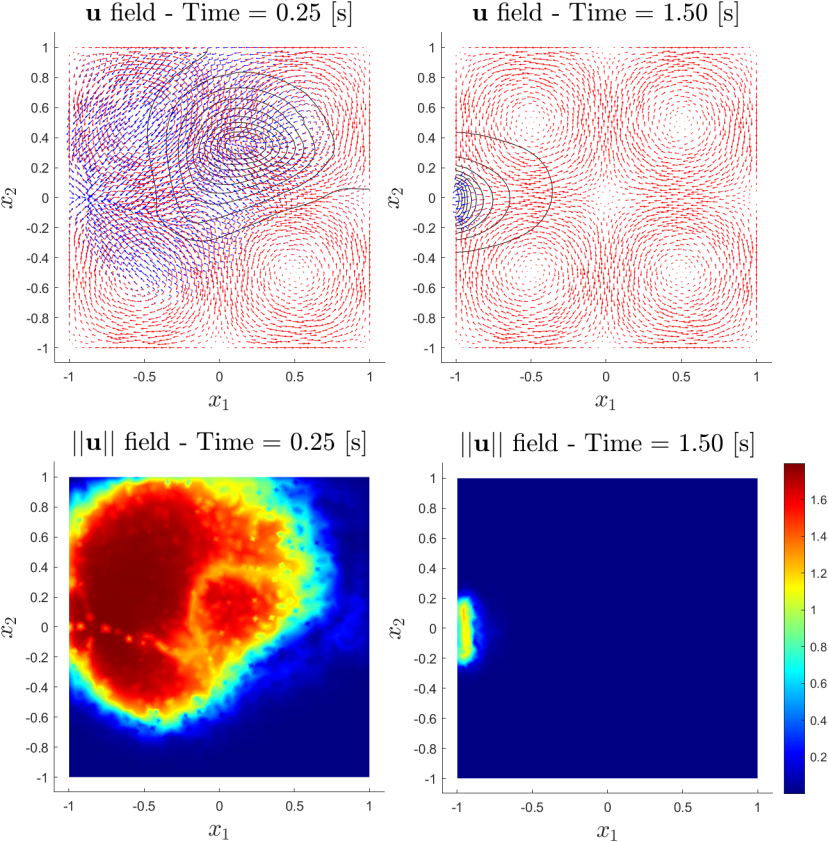}  
\caption{Test Case 1.  
{\it Top:} 
Swarm density contours ({\it black}), 
optimal robot velocity control $\mathbf{u}(\mathbf{x},t)$ ({\it blue}), and fluid velocity field $\mathbf{F}$ ({\it red}) 
at times 
$t=0.25 (s)$ and $t=1.5 (s)$. 
{\it Bottom:} 
Euclidean norm 
$||\mathbf{u}||$ of the optimal control, with a maximum value of $\sqrt{2}$.
}
\label{u_dyn}
\end{figure}

\begin{figure}[h!]
\centering
\includegraphics[width=0.45 \textwidth]{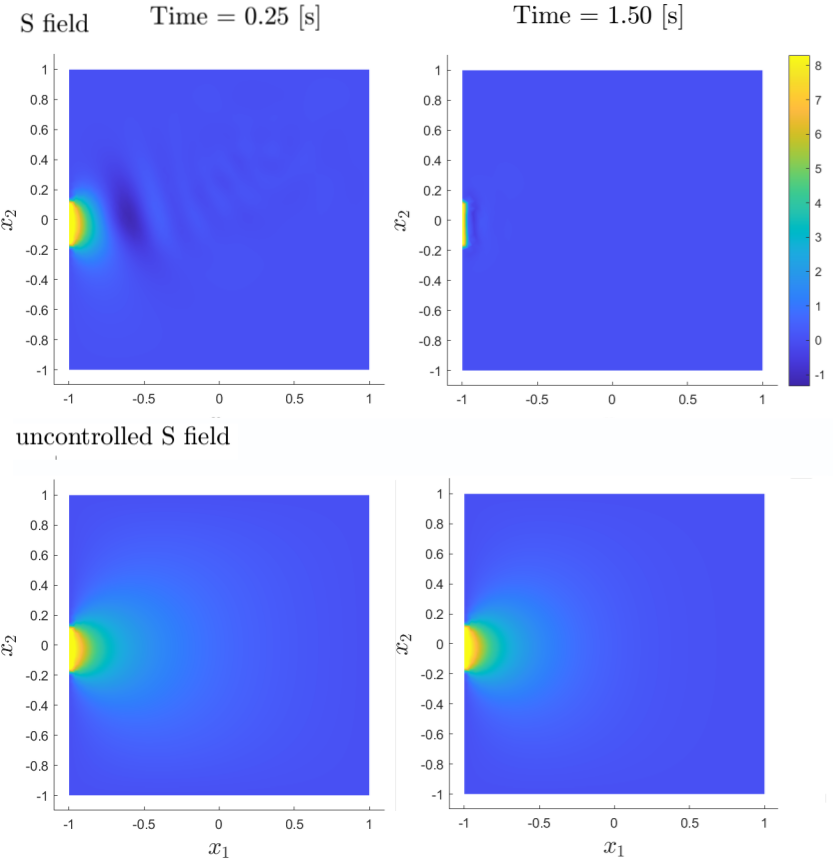}
\caption{Test Case 1. Controlled ({\it top}) and uncontrolled ({\it bottom}) environmental field $S(\mathbf{x},t)$ at times $t=0.25(s)$ and $t=1.5(s)$.   
}
\label{S_dyn}
\end{figure}

\begin{figure}[h!]
\centering
\subfigure{\includegraphics[width=0.45\textwidth]{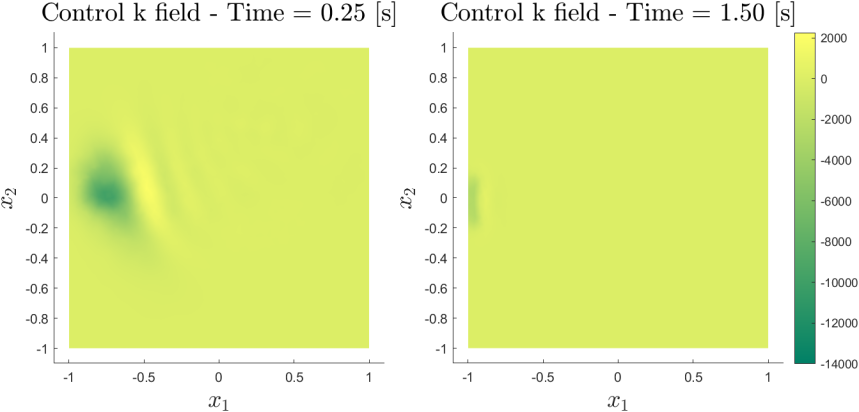}}
\caption{Test Case 1. Optimal actuation $k(\mathbf{x},t)$ at times 
$t=0.25(s)$ and $t=1.5(s)$. The distributed actuation $k$ is correlated with 
the swarm density $q$.}
\label{k_dyn}
\end{figure}

\begin{figure}[h!]
\centering
\subfigure{\includegraphics[width=0.5\textwidth]{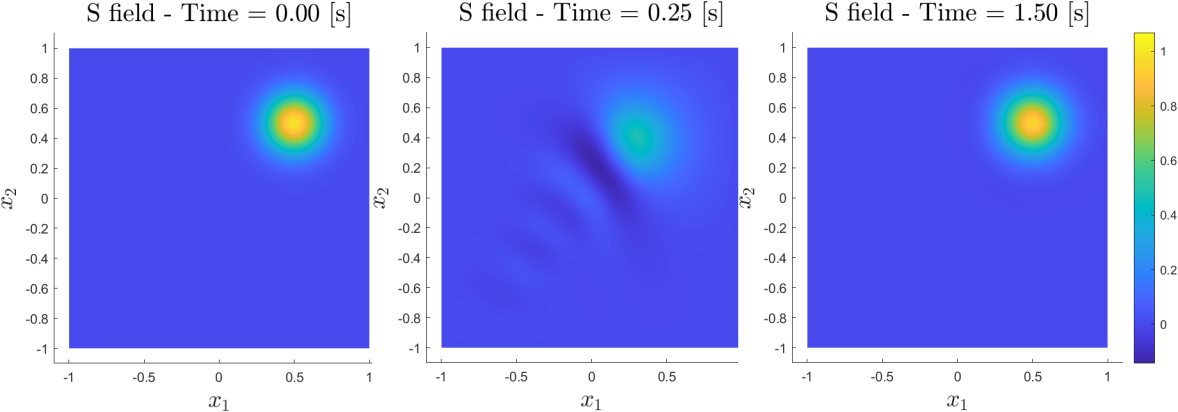}}
\caption{Test Case 2. Controlled environmental field $S(\mathbf{x},t)$ at times $t=0(s)$, $t=0.25(s)$, and $t=1.5(s)$.
For the sake of brevity, the uncontrolled dynamics are not shown since 
they exhibit a trivial diffusion to zero due to the homogeneous Dirichlet boundary condition.}
\label{S_dyn_TC2}
\end{figure}

Two test cases are considered.
In Test Case 1, we solve a regulation problem with target distribution $z=0$ and the Dirichlet boundary condition illustrated in Figure \ref{schem_diag}, with $S = S_d = 10$ along $\Gamma_d$.  
Test Case 2 is 
a tracking problem with $z = S_0(\mathbf{x})$ 
and a homogeneous Dirichlet boundary condition. In Figure \ref{mass_conv}, we compare 
the controlled and 
uncontrolled dynamics of the total mass of the environmental field $S$, defined as $m_S(t) =\int_{\Omega} S(\cdot,t) d\Omega$, for the two test cases. For Test Case 1, the uncontrolled steady-state value of $m_S(t)$ depends on the equilibrium balance between the mass generated along $\Gamma_d$ and the mass absorbed along the rest of the boundary. 
In the controlled case,  however, 
the robotic swarm drives 
$m_S(t)$  to zero through coordinated motion, defined by their  self-propelled velocity $\mathbf{u}$, and the intensity $k$ of their distributed actuation.  In Test Case 2, the uncontrolled mass $m_S(t)$  exponentially converges to zero due to diffusion and the homogeneous Dirichlet boundary condition, while the controlled mass $m_S(t)$ is driven to $\int_{\Omega} S_0(\cdot) d\Omega$ due to the efforts of the swarm to maintain $S$ at its initial condition $S_0(\mathbf{x})$. 
Figure \ref{u_dyn} shows snapshots of the swarm density dynamics under the action of the controls and the fluid velocity 
field for Test Case 1. Figure \ref{S_dyn} compares snapshots of the controlled and uncontrolled dynamics of $S$ for Test Case 1, and Figure \ref{k_dyn} presents snapshots of the optimal actuation $k$ for this case. Finally, snapshots of the controlled dynamics of $S$
for Test Case 2 are shown in Figure \ref{S_dyn_TC2}.


\bibliography{References}               
\end{document}